\documentclass[11pt,a4paper]{amsart}

\usepackage{amsmath,mathrsfs,amssymb}

\newtheorem{theorem}{Theorem}[section]

\newtheorem{proposition}[theorem]{Proposition}
\newtheorem{corollary}[theorem]{Corollary}

\theoremstyle{definition}

\theoremstyle{plain}

\newcommand{\FamilyF}{\mathscr{F}}

\newcommand{\ol}[1]{\overline{#1}}
\newcommand{\loopauto}[2]{\ensuremath{\hat\Gamma_{#1}({#2})}}

\begin{document}
    \title{The Loop Problem for Rees Matrix Semigroups}

\maketitle

\begin{center}
    Mark Kambites \\

    \medskip

    School of Mathematics, \   University of Manchester \\
    Manchester M60 1QD, \  England \\

    \medskip

    \texttt{Mark.Kambites@manchester.ac.uk} \\

\end{center}

\begin{abstract}
We study the relationship between the loop problem of a semigroup,
and that of a Rees matrix construction (with or without zero) over the
semigroup. This allows us to characterize exactly those completely
zero-simple semigroups for which the loop problem is context-free. We
also establish some results concerning loop problems for subsemigroups
and Rees quotients.
\end{abstract}

\bigskip

\section{Introduction}\label{sec_intro}

A key facet of both combinatorial group theory and combinatorial
semigroup theory is the study of the \textit{word problem}, that is,
the problem of deciding algorithmically whether two words in the
generators of a group or semigroup represent the same element \cite{Dehn12}.
In a group, deciding whether two words represent the same element is
algorithmically equivalent to deciding whether one word represents the 
identity. Hence, the word problem can be viewed from an alternative 
perspective, as the \textit{formal language} of all words over the 
generators which represent the identity. Perhaps surprisingly, there are
number of strong results correlating structural properties of groups with
language-theoretic properties of their word problems (see, for example,
\cite{KambitesCounter,Herbst91,Holt05,KambitesWordProblemsRecognisable,Muller83,Stewart99}).

While many major group-theoretic results about the word
problem as a decision problem have been successfully generalised to semigroups
\cite{Belyaev84,Boone74,Murskui67}, and in some
cases beyond \cite{Belegradek96,Kharlampovich95,Kublanovski97,Kukin95}, it is not so clear how one might
apply the language-theoretic approach in a more general setting. In \cite{KambitesLoop}, we introduced and began the study of a new
way of associating to each finitely generated monoid or semigroup a formal
language, called the \textit{loop problem}. In the case of a group the
loop problem is essentially the same as the word problem, and for more
general semigroups, it still encodes the structure of the semigroup in
a similar way.

The \textit{Rees matrix construction}, originally used by Rees \cite{Rees40}
to give a combinatorial description of \textit{completely simple} and \textit{completely zero-simple} semigroups
in terms of their maximal subgroups, has become one of the most pervasive
ideas in semigroup theory. The groups in the construction can be 
replaced with more general semigroups \cite{Ayik99, Descalco01, McAlister83,Petrich04}
or even partial algebras such as small categories and semigroupoids
\cite{KambitesSgpoidPres,KambitesAutoCat,Lawson00}, with the resulting
construction used either to describe an interesting class of semigroups, or
to embed semigroups into other semigroups with nicer properties.

In \cite{KambitesLoop} we started to explore the relationship between the 
loop problem and a special case of the Rees matrix construction, by 
demonstrating a connection between the loop problems of completely simple 
semigroups and the word problems of their maximal subgroups. The main aim 
of the present paper is to develop this study more systematically, by 
considering Rees matrix constructions in more generality. In the process,
we also establish some results regarding the loop problems of subsemigroups
and Rees quotients.

In addition to this introduction, this paper comprises five
sections. In Section~\ref{sec_prelim}, we briefly recall the definition
of the loop problem and the Rees matrix construction, and some basic
results concerning them. Section~\ref{sec_subideal} proves some new
results about the relationship
between the loop problem of a semigroup and that of certain of its
subsemigroups and Rees quotients. Sections~\ref{sec_semitorees}
and~\ref{sec_reestosemi} prove our main theorems about the loop problem
and Rees matrix semigroups. Finally, in Section~\ref{sec_czeros} we consider the
implications of our main results for the class of completely zero-simple
semigroups.

\section{Preliminaries}\label{sec_prelim}

In this section we briefly recall the definition of the loop 
problem of a monoid or semigroup, and of the Rees matrix construction. We 
assume the reader to be familiar with basic definitions and foundational 
results from semigroup theory and from the theory of (both finite and 
infinite) automata and formal languages. A comprehensive introduction to 
the loop problem, including a brief recap of the required 
automata-theoretic prerequisites, can be found in \cite{KambitesLoop}. 
More detailed introductions to automata and formal languages can be found 
in \cite{Hopcroft69,Lawson03,Pin86}, while the standard introductions to 
semigroup theory are \cite{Clifford61,Howie95}.

A \textit{choice of (monoid) generators} for a monoid $M$ is a surjective
morphism $\sigma : X^* \to M$ from a free monoid onto $M$. Similarly, a
\textit{choice of (semigroup) generators} for a semigroup $S$ is a
surjective morphism $\sigma : X^+ \to S$ from a free semigroup onto $S$.
Such a choice of generators is called \textit{finite} if $X$ is finite,
and a semigroup or monoid which admits a finite choice of generators is
called \textit{finitely generated}.

Recall that the \textit{(right) Cayley graph} $\Gamma_\sigma(M)$ of a monoid
$M$ with respect to a choice of generators $\sigma : X^* \to M$ is a directed graph,
possibly with multiple edges and loops, with edges labelled by elements of $X$.
Its vertices are the elements of $M$, and it has an edge from $a \in M$ to
$b \in M$ labelled $x \in X$ exactly if $a (x \sigma) = b$ in the monoid $M$.

Now let $\ol{X} = \lbrace \ol{x} \mid x \in X \rbrace$
be a set of formal inverses for the generators in $X$, and let
$\hat{X} = X \cup \ol{X}$. We extend the map $x \to \ol{x}$ to an
involution on $\hat{X}^*$ by defining $\ol{\ol{x}} = x$ for all $x \in X$,
and $\ol{x_1 \dots x_n} = \ol{x_n} \dots \ol{x_1}$ for all
$x_1, \dots, x_n \in \hat{X}$.
The \textit{(right) loop automaton} $\loopauto{\sigma}{M}$ of $M$ with
respect to $X$ is obtained from the Cayley graph $\Gamma_\sigma(M)$ by
adding for each edge labelled $x$ an inverse edge, in the opposite
direction, labelled $\ol{x}$.
Notice that in the loop automaton,
$w \in \hat{X}$ labels a path from $p$ to $q$ if and only if $\ol{w}$
labels a path from $q$ to $p$. We view the loop automaton as a (typically
infinite) automaton over $\hat{X}$, with start state and terminal
state the identity of $M$.
The \textit{(right) loop problem} of $M$ with respect to $\sigma$ is the
language $L_\sigma(M) \subseteq \hat{X}^*$ of words recognised by
the loop automaton $\loopauto{\sigma}{M}$.

Words in $X^*$ and $\ol{X}^*$ are called \textit{positive} and
\textit{negative} words respectively; words in $\hat{X}^*$ which are
neither positive nor negative are called \textit{mixed}. Similarly, an
edge or path in $\loopauto{\sigma}{M}$ is called \textit{positive}
[\textit{negative}, \textit{mixed}] if it has a positive [respectively
negative, mixed] label.

Now let $S$ be a semigroup, and let $S^1$ be the monoid obtained from
$S$ by adjoining an identity $1$ (even if $S$ already has an identity).
Then a choice of semigroup generators $\sigma : X^+ \to S$ extends
uniquely to a choice of monoid generators
$\sigma^1 : X^* \to S^1$. We define the \textit{loop automaton}
$\loopauto{\sigma}{S}$ and the \textit{loop problem} $L_\sigma(S)$
of $S$ with respect to $\sigma$ to be respectively the loop automaton
$\loopauto{\sigma^1}{S^1}$ and the loop problem $L_{\sigma^1}(S^1)$.

We now turn our attention to the Rees matrix construction. Let $S$ be
a semigroup, $I$ and $J$ sets, $0$ a symbol not in $S$ and $P$ an
$J \times I$ matrix with entries drawn from $S \cup \lbrace 0 \rbrace$.
The \textit{Rees matrix semigroup with zero}
$$M^0(S; I, J; P)$$ is the semigroup with set of elements
$$(I \times S \times J) \cup \lbrace 0 \rbrace$$
and multiplication given by
$$(i_1, g_1, j_1) (i_2, g_2, j_2) =
   \begin{cases}
      (i_1, g_1 P_{j_1 i_2} g_2, j_2) &\text{ if } P_{j_1 i_2} \in S \\
      0 &\text{ if } P_{j_1 i_2} = 0.
   \end{cases}
$$
If $P$ contains no $0$ entries, $I \times S \times J$ forms a subsemigroup
of $M^0(S; I, J; P)$, called a \textit{Rees matrix semigroup (without zero)}
and denoted $M(S; I,J; P)$.

\section{Subsemigroups and Ideals}\label{sec_subideal}

In this section we prove some basic results relating the loop problem of
a semigroup to the loop problems of certain of its subsemigroups and quotients.

Let $S$ be a semigroup with a subsemigroup $T$ such that $STS \subseteq T$;
such a subsemigroup $T$ is called an \textit{ideal} of $S$. The \textit{Rees
quotient} of $S$ by $T$ is the semigroup with set of elements
$S \setminus T \cup \lbrace 0 \rbrace$ and multiplication given by
$$st = \begin{cases}
    0 &\text{ if } s = 0 \text{ or } t = 0; \\
    \text{the $S$-product $st$} &\text{ otherwise}.
\end{cases}
$$
Some authors use the symbol $T$ to denote the zero element $0$.
The map 
$$S \to S/T, \  s \mapsto 
\begin{cases}
    s &\text{ if } s \notin T \\
    0 &\text{ if } s \in T.
\end{cases}$$
is clearly a morphism from $S$ onto $S/T$.
If $\sigma : X^+ \to S$ is a choice of generators for $S$, then composition
with this map induces a choice of generators for $S/T$, which we denote
by $\sigma / T : X^+ \to S/T$. The map $S \to S / T$ also induces a
label-preserving graph morphism from $\loopauto{\sigma}{S}$ onto
$\loopauto{\sigma/T}{S/T}$.

Recall that a subset of a semigroup $S$ is called \textit{rational} if $S$ 
is the image in $S$ of a regular language under some morphism from a 
finitely generated free monoid to $S$, or equivalently, under every 
morphism from a finitely generated free monoid to $S$. A rational subset 
which is also an ideal is called a \textit{rational ideal}.

\begin{proposition}\label{prop_reesquotient}
Let $\FamilyF$ be a family of languages closed under product, union,
Kleene closure and division by finite [respectively, regular] languages.
Suppose $\sigma : X^* \to S$ is a choice of generators for a
semigroup $S$, and that $L_\sigma(S)$ belongs to $\FamilyF$. If $T$
is a finite [respectively, rational] ideal of $S$ then $L_{\sigma/T}(S/T)$
belongs to $\FamilyF$.
\end{proposition}
\begin{proof}
Clearly, a path
from $1$ to $1$ in $\loopauto{\sigma/T}{S/T}$ consists of either
\begin{itemize}
\item[(i)] a path from $1$ to $1$ which does not visit $0$, and hence which
is the image of a path from $1$ to $1$ in $\loopauto{\sigma}{S}$ with the same label; or
\item[(ii)] the concatenation of
\begin{itemize}
\item a path from $1$ to $0$, which lifts to a path with the same label from $1$ to an element
of $T$ in $\loopauto{\sigma}{S}$,
\item zero or more paths from $0$ to $0$, each of which lifts to a path from
an element of $T$ to an element of $T$; and
\item a path from $0$ to $1$, which lifts to a path with the same label
from an element of $T$ to $1$ in $\loopauto{\sigma}{S}$.
\end{itemize}
\end{itemize}
It follows that
$$L_{\sigma/T}(S/T) \ = \ L_\sigma(S) \cup L_{1T} L_{TT} L_{T1}$$
where $L_{1T}$, $L_{TT}$ and $L_{T1}$ denote the sets of words
which label paths in $\loopauto{\sigma}{S}$ from $1$ to an element of
$T$, between two elements of $T$, and from an element of $T$ to $1$,
respectively. 

Since $T$ is finite [rational], there exists a finite [respectively,
regular] language $R \subseteq X^+$ with $R \sigma = T$. Let
$\ol{R} = \lbrace \ol{w} \mid w \in R \rbrace$; then $\ol{R}$ is also
finite [respectively, regular].
Notice that, for any word $r \in X^+$, we have that $r$ labels a path from
$1$ to an element of $T$ exactly if $r \in R$, and hence that $\ol{r}$
labels a path from an element of $T$ to $1$ if and only $\ol{r} \in \ol{R}$.
It follows that
$$L_{1T} = \lbrace w \in \hat{X}^* \mid w\ol{r} \in L_\sigma(S) \text{ for some } r \in R \rbrace = L_\sigma(S) \ol{R}^{-1}.$$
Similarly, we have
and
$$L_{T1} = \lbrace w \in \hat{X}^* \mid r w \in L_\sigma(S) \text{ for some } r \in R \rbrace = R^{-1} L_\sigma(S)$$
and
$$L_{TT} = \lbrace w \in \hat{X}^* \mid r_1 w \ol{r_2} \in L_\sigma(S) \text{ for some } r_1, r_2 \in R \rbrace = R^{-1} L_\sigma(S) \ol{R}^{-1}.$$
The claim now follows immediately from the presumed closure properties of
$\mathscr{F}$.
\end{proof}

Next, we consider loop problems of subsemigroups.
We say that a subsemigroup $T$ of $S$ is \textit{pseudo-right-unitary} if for
every element $a \in S$ there exists an element $b \in T$ such that
for each element $x \in T$ with $ax \in T$, we have $ax = bx$.
Intuitively, $T$ is pseudo-right-unitary if the partial action of $S$
on $T$ by left multiplication is no more complex than the action of
$T$ on itself by left multiplication. Note in particular that if $T$ is
right unitary then $x \in T$ and $ax \in T$ together imply that $a \in T$,
so that right unitary subsemigroups (and in particular unitary subsemigroups)
are pseudo-right-unitary.

In fact, a weaker condition will suffice for our purpose. A subsemigroup
is \textit{weakly pseudo-right-unitary} if for every element $a \in S$ and
pair of elements $x,y \in T$ such that $ax \in T$, there exists $b \in T$
with $ax = bx$ and $ay = by$.

\begin{proposition}\label{prop_pru}
Let $T$ be a finitely generated weakly pseudo-right-unitary subsemigroup of
a semigroup $S$. Let $\sigma : X^+ \to T$ and $\tau : Y^+ \to S$ be
choices of generators such that $X \subseteq Y$ and $\sigma$ is the
restriction of $\tau$ to $X^*$. Then $L_\sigma(T) = L_\tau(S) \cap \hat{X}^+$.
\end{proposition}
\begin{proof}
The proof is similar to that of \cite[Theorem 5.4]{KambitesLoop}.
Clearly, $\loopauto{\sigma}{T}$ is embedded into 
$\loopauto{\tau}{S}$, so that any word accepted by the former is
also accepted by the latter, and $L_\sigma(T) \subseteq L_\tau(S) \cap \hat{X}^*$.

Conversely, suppose $w$ is a word in $L_\tau(S) \cap \hat{X}^*$. Then
$\loopauto{\tau}{S} = \loopauto{\tau^1}{S^1}$ has a loop at $1$ labelled
$w$. Write $w = u_0 \ol{v_1} u_1 \ol{v_2} \dots u_{n-1} \ol{v_n}$ with each
$u_i, v_i \in X^*$. By the Zig Zag Lemma \cite[Lemma~4.1]{KambitesLoop},
there exist $p_0, \dots, p_n \in S^1$ such that $p_0 = p_n = 1$ and
$p_i (u_i \tau) = p_{i+1} (v_{i+1} \tau)$ for
$0 \leq i < n$.

We claim that the elements $p_0, \dots, p_n$ can all be
chosen to lie in $T \cup \lbrace 1 \rbrace$. Indeed, suppose we are give
$p_0, \dots, p_n$ satisfying the above equations and not all lying in
$T \cup \lbrace 1 \rbrace$, and let $k$ be
minimal such that $p_k \notin T \cup \lbrace 1 \rbrace$. Certainly,
$1 \leq k \leq n-1$.
Now we have $u_{k-1} \tau \in T$ and $p_{k-1} \in T \cup \lbrace 1 \rbrace$,
so that $p_{k-1} (u_{k-1} \tau) = p_k (v_k \tau)$ lies in $T$. Since
$v_k \tau$ and $u_k \tau$ also lie in $T$ and $T$ is weakly pseudo-right-unitary,
there exists an element $q_k \in T$ with
$q_k (v_k \tau) = p_k (v_k \tau) = p_{k-1} (u_{k-1} \tau)$ and
$q_k (u_k \tau) = p_k (u_k \tau) = p_{k+1} (v_{k+1} \tau)$. Replacing
$p_k$ with $q_k$, we obtain a sequence $p_0, \dots, p_n$ with
strictly fewer elements outside $T \cup \lbrace 1 \rbrace$, and continuing
in this way we eventually obtain a sequence with the desired properties.
Another application of the Zig Zag Lemma shows that
$\loopauto{\sigma}{T} = \loopauto{\sigma^1}{T^1}$ has a
loop at $1$ labelled $w$, so that $w \in L_\sigma(T)$, as required.
\end{proof}

We denote by $S^0$ the semigroup obtained from $S$ by adjoining a zero.
Since $S$ is a right unitary subsemigroup of $S^0$, Proposition~\ref{prop_pru}
has the following immediate corollary.

\begin{corollary}\label{cor_removezero}
Let $\sigma : X^+ \to S$ be a choice of generators for a semigroup $S$.
Let $Y = X \cup \lbrace z \rbrace$ and let $\tau : Y^+ \to S^0$ be the
(unique) extension of $\sigma$ to a choice of generators for $S^0$. Then
$L_\sigma(S) = L_\tau(S^0) \cap \hat{X}^*$.
\end{corollary}

Conversely, we can also obtain the loop problem of $S^0$ from that of
$S$, but to do so we need some more operations.

\begin{theorem}\label{thm_adjoinzero}
Let $\FamilyF$ be a family of languages closed under union, product, Kleene
closure and rational transductions, and let $S$ be a semigroup. Then
the loop problem for $S$ belongs to $\FamilyF$ if and only if the loop
problem for $S^0$ belongs to $\FamilyF$.
\end{theorem}
\begin{proof}
Let $\sigma : X^+ \to S$ be a choice of generators for $S$, let
$Y = X \cup \lbrace z \rbrace$ and let $\tau : Y^+ \to S^0$ be the
unique extension of $\sigma$ to a choice of generators for $S^0$.
Since intersection with a regular language can be performed by a
rational transduction, one implication is immediate from
Corollary~\ref{cor_removezero}

Conversely, the loop automaton $\loopauto{\tau}{S^0}$ consists of the
loop automaton $\loopauto{\sigma}{S}$ with an extra vertex $0$ adjoined.
In addition to the edges in $\loopauto{\sigma}{S}$, it has an edge from
every vertex to $0$ labelled $z$, a corresponding edge from $0$ to every
vertex labelled $\ol{z}$, and an edge from $0$ to $0$ labelled $x$ for
each $x \in \hat{X}$.

We define languages
$$L_{10} \ = \ \lbrace u z \mid u \text{ is a prefix of a word in } L_\sigma(S) \rbrace,$$
$$L_{00} \ = \ \lbrace \ol{z} u z \mid u \text{ is a factor of a word in } L_\sigma(S) \rbrace, \text{ and}$$
$$L_{01} \ = \ \lbrace \ol{z} u \mid u \text{ is a suffix of a word in } L_\sigma(S) \rbrace.$$
It is a routine exercise to verify that $L_{10}$, $L_{00}$ and $L_{01}$ can
be obtained using rational transductions from $L_\sigma(S)$, and hence lie in
$\FamilyF$. One can also check that
\begin{itemize}
\item $L_{10}$ is the set of all words labelling paths in $\loopauto{\tau}{S^0}$
      from $1$ to $0$ which do not visit $0$ except at the end;
\item $L_{00}$ is the set of all words labelling paths in $\loopauto{\tau}{S^0}$
      from $0$ to $0$ which do not visit $0$ except at the beginning and end; and
\item $L_{01}$ is the set of all words labelling paths in $\loopauto{\tau}{S^0}$
      from $0$ to $1$ which do not visit $0$ except at the beginning.
\end{itemize}

Now if $\pi$ is a path from $1$ to $1$ in $\loopauto{\tau}{S^0}$, then
$\pi$ either does or does not visit the vertex $0$. If $\pi$ does not
visit zero then it is entirely contained in $\loopauto{\sigma}{S}$ and
its label lies in $L_\sigma(S)$. If $\pi$ does visit $0$ then it must
decompose into a path from $1$ to $0$ which does not visit $0$
except at the end, followed by zero or more paths from $0$ to $0$ which do
not visit $0$ except at the end, followed by a path from $0$ to $1$ which
does not visit $0$. It follows that
$$L_\tau(S^0) \ = \ L_\sigma(S) \cup L_{10} L_{00}^* L_{01}.$$
Thus, $L_\tau(S^0)$ can be obtained from languages in $\FamilyF$ using
the operations of union, product and Kleene closure, and so itself lies
in $\FamilyF$.
\end{proof}

\section{From Semigroups to Rees Matrix Semigroup}\label{sec_semitorees}

In this section, we show how to describe the loop problem of a finitely
generated Rees matrix semigroup in terms of the loop problem of the
underlying semigroup. We first do this for Rees matrix constructions
without zero; we shall subsequently apply one of the results of
Section~\ref{sec_subideal} extend this result to the case of Rees matrix
constructions with zero.

\begin{theorem}\label{thm_semitorees}
Let $S$ be a semigroup and $M = M(S; I, J; P)$ a finitely generated 
Rees matrix semigroup without zero. Then the loop problem of $M$ is the 
Kleene closure of a rational transduction of the loop problem of $S$. 
\end{theorem}
\begin{proof}
By \cite[Main Theorem]{Ayik99}, we may assume that $I$ and $J$ 
are finite and $S$ is finitely generated. Let $\sigma : X^+ \to S$ be a 
finite choice of generators for $S$, and $\tau : Y^+ \to M$ a finite 
choice of generators for $M$. We shall show that $L_\tau(M)$ 
is a rational transduction of $L_\sigma(S) = L_{\sigma^1}(S^1)$; the
argument is a refinement of that used to prove
\cite[Theorem~5.5]{KambitesLoop}.

For each $y \in Y$, suppose $y \tau = (i_y, g_y, j_y)$ and let
$w_y \in X^+$ be a word representing $g_y \in S$. For each
$i \in I$ and $j \in J$, let $w_{ji} \in X^+$ be a word
representing $P_{ji} \in S$.

We define a finite state transducer from $\hat{X}^*$ to $\hat{Y}^*$ with
\begin{itemize}
\item vertex set $(I \times J) \ \cup \ \lbrace A, Z \rbrace$ where
$A$ and $Z$ are new symbols;
\item initial state $A$;
\item terminal state $Z$;
\item for each generator $y \in Y$, an edge from
$A$ to $(i_y,j_y)$ labelled $(w_y, y)$;
\item for each generator $y \in Y$, an edge from
$(i_y,j_y)$ to $Z$ labelled $(\ol{w_y}, \ol{y})$;
\item for each generator $y \in Y$, each $k \in J$ and each $i \in I$,
an edge from $(i,k)$ to $(i,j_y)$ labelled
$(w_{ki_y} w_y, y)$; and
\item for each generator $y \in Y$, each
$j \in J$ and each $i \in I$, an edge from
$(i,j_y)$ to $(i,j)$ labelled $(\ol{w_y} \ \ol{w_{ji_y}}, \ol{y})$.
\end{itemize}

We say that a path starting at $1$ in $\loopauto{\tau}{M}$ is
\textit{non-returning} if it does not visit the vertex $1$ except
at the start and (possibly) the end.
Now let $g \in S$, $i \in I$, $j \in J$ and $v \in \hat{Y}^+$ and suppose
$n$ is a positive integer. We claim that the following conditions are
equivalent:
\begin{itemize}
\item[(i)] the loop automaton $\loopauto{\tau}{M}$ has a non-returning path
           of length $n$ from $1$ to $(i, g, j)$ labelled $v$;

\item[(ii)] the transducer has a path of length $n$ from $A$ to $(i,j)$
            labelled $(u, v)$ for some $u \in \hat{X}^+$ such that the
            loop automaton $\loopauto{\sigma}{S}$ has a non-returning
            path of length $n$ from $1$ to $g$ labelled $u$.
\end{itemize}
We prove this claim by induction on $n$. The case $n=1$ is immediate from
the definition of the transducer, so suppose $n > 1$ and that the claim holds
for smaller $n$.

Suppose first that (i) holds, and let $\pi$ be the path given by the
hypothesis. Let $e$ be the last edge the path $\pi$, and let $\pi'$ be
the path $\pi$ with the last edge removed, so that $\pi = \pi' e$. Let
$v'$ be the label of $\pi'$. Since $n > 1$ and $\pi$ is non-returning,
the path $\pi'$ must end at a vertex of the form $(i',g',k)$.
It follows easily from the definition of
the multiplication in a Rees matrix semigroup that the vertices in the
loop automaton corresponding to elements with first coordinate $i$ are
connected to the rest of the automaton only via the vertex $1$.
Hence, since the path $\pi$ is non-returning, we
must have $i = i'$, so that $\pi$ actually ends at $(i,g',k)$. Now $\pi$
is a path of length $n-1$, so by
the inductive hypothesis, the transducer has a path of length $n-1$ from
$A$ to $(i,k)$ labelled $(u', v')$ for some word $u' \in X^+$ such that
$\loopauto{\sigma}{S}$ has a path from $1$ to $g'$ labelled $u'$.

We treat separately the case where $e$ is a positive edge, and that
where $e$ is a negative edge.
Suppose first that $e$ is a positive edge,
with label $y \in Y$. Then clearly $j = j_y$ and by definition the
transducer
has an edge from $(i, k)$ to $(i,j)$ with label $(w_{ki_y} w_y, y)$. Hence,
the transducer has a path of length $n$ from $1$ to $(i,j)$ with label
$$(u', v') (w_{ki_y} w_y, y) \ = \ (u' w_{ki_y} w_y, v' y) \ =
 \ (u' w_{ki_y} w_y, v).$$
Let $u = u' w_{ki_y} w_y$. Now from the definition of the loop automaton,
we must have
$$(i, g, j) \ = \ (i,g',k) (y \tau) \ = \ (i, g', k) (i_y, g_y, j_y) \ =
 \ (i, g' P_{k i_y} g_y, j_y).$$
Equating second coordinates, we see that $g = g' P_{ki_y} g_y$.
Now since $w_{ki_y} w_y \in X^+$ represents $P_{ki_y} g_y$, it follows that
$\loopauto{\sigma}{S}$ has a path from $g'$ to $g$ labelled $w_{ki_y} w_y$.
Combining with the path given by the inductive hypothesis, we see that
$\loopauto{\sigma}{S}$ has a path from $1$ to $g$ labelled
$u' w_{ki_y} w_y = u$, which shows that (ii) holds in the case that
$e$ is a positive edge.

Suppose now that $e$ is a
negative edge, with label $\ol{y}$ for some $y \in Y$. In this case
$k = j_y$ and the transducer has an edge from
$(i,k)$ to $(i,j)$ with label
$(\ol{w_y} \ \ol{w_{ji_y}}, \ol{y})$. Hence, the the transducer has a
path of length $n$ from $A$ to $(i,j)$ with label
$$(u', v') (\ol{w_y} \ \ol{w_{ji_y}}, \ol{y}) \ = \ (u' \ \ol{w_y} \ \ol{w_{ji_y}}, v' \ \ol{y}) \
 = \ (u' \ \ol{w_y} \ \ol{w_{ji_y}}, v).$$
Let $u = u' \ol{w_y} \ \ol{w_{ji_y}}$.
Now from the definition of the loop automaton, we must have
$$(i, g', k) \ = \ (i,g',j) (y \tau) \ = \ (i, g', j) (i_y, g_y, j_y)
 \ = \ (i, g' P_{j i_y} g_y, j_y).$$
Again equating second coordinates, we see this time that
$g P_{ji_y} g_y = g'$. Since $w_{j i_y} w_y$ represents $P_{ji_y} g_y$, it
follows this time that $L_\sigma(S)$ has a path from $g$ to $g'$ labelled
$w_{ji_y} w_y$, and hence a path from $g'$ to $g$ with label
$$\ol{w_{ji_y} w_y} \ = \ \ol{w_y} \ \ol{w_{ji_y}}.$$
Combining with the path given by the inductive hypothesis, we obtain a
path from $1$ to $g$ labelled $u = u' \ol{w_y} \ \ol{w_{ji_y}}$, so that
(ii) again holds.

Conversely, suppose (ii) holds, and let $\pi$ be a path of length
$n$ in the transducer from $A$ to $(i,j)$ labelled $(u,v)$ for some $u \in X^+$
such that $\loopauto{\sigma}{S}$ has a path from $1$ to $g$ with label $u$.
Let $e$ be the last edge of $\pi$ and let
$\pi'$ be the path $\pi$ with the final edge removed. Then $\pi'$ is
a path of length $n-1$ from $A$ to some vertex $(i',k)$ with label of the
form $(u', v')$. Moreover, it follows easily from the definition of the
transducer that $i = i'$, so that $\pi'$ ends at $(i, k)$.
Let $g' \in G$ be the element represented by $u'$. Then
by the inductive hypothesis,
there exists a path of length $n-1$ in the loop automaton from $1$ to
$(i,g',k)$ with label $v'$. Now $e$ is an edge from $(i,k)$ to
$(i,j)$. From the definition of the edges in the transducer, we see that there
exists $y \in Y$ such that either $j_y = j$ and $e$ has
label $(w_{ki_y} w_y, y)$, or else $j_y = k$ and $e$ has
label $(\ol{w_y} \ \ol{w_{ji_y}}, \ol{y})$. As before,
we treat these two cases separately.

In the former case, observe that we have $u = u' w_{ki_y} w_y$
from which we deduce that $g = g' P_{ki_y} g_y$. But now
$$(i,g',k) (y \tau) = (i, g' P_{ij_y} g_y, j_y) = (i, g, j)$$
so we see that the loop automaton has an edge from $(i,g',k)$ to
$(i,g,j)$ labelled $y$. Combining this with the path whose existence
we deduced using
the inductive hypothesis, we conclude that the loop automaton has a path
from $1$ to $(i,g,j)$ labelled $v = v' y$, so that (i) holds as required.

Next we consider the case in which $j_y = k$ and $e$ has label
of the form  $(\ol{w_y} \ \ol{w_{ji_y}}, \ol{y})$.
Here we have $u = u' \ \ol{w_y} \ \ol{w_{ji_y}}$, so $\loopauto{\sigma}{S}$
has a path from $g'$ to $g$ with label $\ol{w_y} \ \ol{w_{ji_y}}$, and hence
a path from $g$ to $g'$ with label $w_{ji_y} w_y$. It follows that
$$(i,g,k) (y \tau) = (i, g P_{ji_y} g_y, j_y) = (i, g', j)$$
so we see that the $\loopauto{\tau}{M}$ has an edge from $(i,g,j)$ to
$(i,g',k)$ labelled $y$, and hence an inverse edge from $(i,g',k)$ to
$(i, g, j)$ labelled $\ol{y}$. Combining this with the path whose existence we deduced
using the inductive hypothesis, we conclude that the loop automaton has a 
path from $1$ to $(i,g,j)$ labelled $v = v' \ol{y}$, so that (i) holds as required.
This completes the proof that (i) and (ii) are equivalent.

Now let $K$ be the set of all non-empty words labelling non-returning loops at
$1$ in $\loopauto{\tau}{M}$. We claim that $K$ is exactly the image of
$L_\sigma(S)$ under the transduction defined by our transducer. Since
$L_\tau(M)$ is clearly the Kleene closure of $K$, this will suffice to
complete the proof.

Suppose first that $v \in K$. Then by definition the loop automaton
$\loopauto{\tau}{M}$ has a non-returning loop at $1$ labelled $v$. Note that
all edges in $\loopauto{\tau}{M}$ which end at $1$ run from vertices
corresponding to generators
$y$ and have label $\ol{y}$, so we may assume that the last edge of the path
runs from $y \tau$ to $1$, and has label $\ol{y}$. Let $\pi$ be
the path
without this last edge, so that $\pi$ runs from $1$ to $y \tau$, and
let $v' \in \hat{Y}^+$ be the label of this path, so that $v = v' \ol{y}$.
Then by the equivalence
above, the transducer has a path from $A$ to $(i_y, j_y)$ labelled $(u,v')$
for some $u \in \hat{X}^+$ such that $u$ labels a path in
$\loopauto{\sigma}{S}$ from $1$ to $g_y$. But directly from the definition,
the transducer also has an edge from $(i_y, j_y)$ to $Z$ labelled
$(\ol{w_y}, \ol{y})$.
Hence, we deduce that $(u \ol{w_y}, v)$ is accepted by the transducer,
where $u \ol{w_y}$ labels a path from $1$ to $1$ in $\loopauto{\sigma}{S}$.
Thus, $v$ lies in the image under the transduction of $L_{\sigma}(S)$.

Conversely, suppose $u \in \hat{X}^+$ lies in the loop problem
$L_{\sigma}(S)$, and that the transducer accepts $(u,v)$. Then the
transducer has a path $\pi$ from $A$
to $Z$ labelled $(u,v)$. Again, we proceed by letting $\pi'$ be the path
obtained from $\pi$ by deleting the last edge. Then there must exist a
generator $y$ such that $\pi'$ ends at $(i_y,j_y)$. Moreover, $\pi'$ must
be labelled $(u', v')$ where $u = u' \ol{w_y}$ and $v = v' \ol{y}$. Now
since $u \in L_{\sigma}(S)$ and $w_y$ represents $g_y$, we deduce that
$u'$ labels a path from $1$ to $g_y$ in $L_\sigma(S)$. By the equivalence above, it follows that the
loop automaton $\loopauto{\tau}{S}$ has a non-returning path from $1$ to
$(i_y, g_y, j_y)$ labelled $v'$. Now it certainly also
has an edge from $(i_y, g_y, j_y)$ to $1$ labelled $\ol{y}$, so we deduce
that $v = v' \ol{y} \in K$, which completes the proof.
\end{proof}

We now use some results from Section~\ref{sec_subideal} to establish a
result corresponding to Theorem~\ref{thm_semitorees} for the case of
Rees matrix constructions with zero.
\begin{theorem}\label{thm_semitoreeszero}
Let $\FamilyF$ be a family of languages closed under rational transduction,
union, product and Kleene closure. Let $S$ be a semigroup and let
$M = M^0(S; I, J; P)$ be a finitely generated Rees matrix semigroup with
zero over $S$. If the loop problem for $S$ belongs to $\FamilyF$ then the
loop problem for $M$ belongs to $\FamilyF$.
\end{theorem}
\begin{proof}
Let $S^0$ denote the semigroup $S$ with an additional zero element $0$
adjoined, and let $M'$ be the Rees matrix semigroup (without zero)
$M' = M(S^0; I, J; P)$. By Theorem~\ref{thm_adjoinzero}, the loop problem for
$S^0$ belongs to $\FamilyF$, and so by Theorem~\ref{thm_semitorees}, the loop
problem for $M'$ also belongs to $\FamilyF$.

Now let
$$T \ = \ I \times \lbrace 0 \rbrace \times J \ \subseteq \ M'.$$
Since $M$ is finitely generated, we deduce by \cite[Main Theorem]{Ayik99} that $I$ and
$J$ are finite, and hence that $T$ is finite. It is readily verified that $T$ is
an ideal of $M'$ and that $M$ is isomorphic to the Rees quotient $M' / T$.

Now since the operations of union and division by finite languages can
easily be realised as rational transductions, we may assume that
$\FamilyF$ is closed under these operations as well as product and
Kleene closure. Thus, Proposition~\ref{prop_reesquotient} tells us that
the loop problem for $M' / T$, and hence that for $M$, belongs to $\FamilyF$
as required.
\end{proof}

\section{From Rees Matrix Semigroups to Subsemigroups}\label{sec_reestosemi}

Our objective in this section is to describe the loop problem for a
semigroup $S$ in terms of the loop problem of a Rees matrix semigroup
$M = M(S; I, J; P)$ or $M = M^0(S; I, J; P)$. In the most general case this
is not possible, since $M$ may not contain sufficient information about $S$.
For example, if the entries in the sandwich matrix $P$ are all drawn from
an ideal of $S$, then $M$ will clearly retain no information about the
structure of $S$ outside this ideal, beyond its cardinality. Thus, for our
results in this section it will be necessary to impose some restrictions
on the elements of the sandwich matrix, and the ideals they generate.

\begin{proposition}\label{prop_reespru}
Let $M = M(S; I, J; P)$ or $M = M^0(S; I, J; P)$ be a finitely generated
Rees matrix semigroup. Choose $i \in I$ and $j \in J$ with $P_{ji} \neq 0$
and let
$$T \ = \ \lbrace (i, s, j) \mid s \in S \rbrace \ \subseteq \ M.$$
If for every $j' \in J$ we have either $P_{j'i} = 0$ or
$S P_{j'i} \subseteq S P_{ji}$ then $T$ is a pseudo-right-unitary
subsemigroup of $M$.
\end{proposition}
\begin{proof}
Clearly, $T$ is a subsemigroup of $M$. Now let $a \in M$. To show that
$T$ is pseudo-right-unitary, we must find an element $b \in T$ such that
for every
element $x \in T$ with $ax \in T$, we have $ax = bx$. If $M$ is a Rees
matrix semigroup with zero and $a = 0$ then $ax$ can never lie in 
$T$, so any element $b \in T$ will fulfil the condition vacuously. So
assume $a \neq 0$, say $a = (i_a, s_a, j_a)$. 

Now by assumption, either $P_{j_a i} = 0$, or $S P_{j_a i} \subseteq S P_{ji}$.
In the former case we have $ax = 0$ for all $x \in T$, so that choosing
any $b \in T$ again suffices vacuously to fulfil the condition.

In the latter case, where $S P_{j_a i} \subseteq S P_{ji}$, we can choose
an element $s_b \in S$ such that $s_a P_{j_a i} = s_b P_{ji}$. Let $b = (i, s_b, j)$. Now suppose
$x = (i, s_x, j) \in T$ is such that
$$ax \ = \ (i_a, s_a P_{j_a i} s_x, j) \ \in \ T.$$
Then from the definition of $T$, we have $i_a = i$. Now by the choice
of $s_b$, we have
$$bx \ = \ (i, s_b P_{ji} s_x, j) \ = \ (i, s_a P_{j_ai} s_x, j) \ = \ ax$$
as required to complete the proof.
\end{proof}

Combining Propositions~\ref{prop_pru} and~\ref{prop_reespru} immediately
yields the following.
\begin{theorem}\label{thm_reestosemi}
Let $M = M(S; I, J; P)$ or $M = M^0(S; I, J; P)$ be a finitely generated
Rees matrix semigroup. Choose $i \in I$ and $j \in J$ with $P_{ji} \neq 0$
and such that for every
$j' \in J$ we have $P_{j'i} = 0$ or $S P_{j'i} \subseteq S P_{ji}$ and let
$$T \ = \ \lbrace (i, s, j) \mid s \in S \rbrace \ \subseteq \ M.$$
Let $\sigma : X^+ \to T$ and $\tau : Y^+ \to M$ be
choices of generators such that $X \subseteq Y$ and $\sigma$ is the
restriction of $\tau$ to $X^+$. Then $L_\sigma(T) = L_\tau(M) \cap \hat{X}^*$.
\end{theorem}

Suppose now that $S$ is a monoid with identity $1$. Recall that an
element $g \in S$ is
called a \textit{unit} if there exists an element $g^{-1} \in S$ with
$gg^{-1} = g^{-1} g = 1$; the collection of units in $S$ forms a subgroup,
called the \textit{group of units} of $S$.
As a corollary of Theorem~\ref{thm_reestosemi} we obtain the following.

\begin{theorem}\label{thm_reestosubsemi}
Let $M = M(S; I, J; P)$ or $M = M^0(S; I, J; P)$ be a finitely generated
Rees matrix semigroup over a monoid $S$. If the sandwich matrix $P$ contains
a unit then there exist choices of generators for $S$ and $M$ such that
the loop problem for $S$ is the intersection of the loop problem for $M$
with a regular language.
\end{theorem}
\begin{proof}
Choose $i \in I$ and $j \in J$ such that $P_{ji}$ is a unit, and let
$$T = \lbrace (i, s, j) \mid s \in S \rbrace.$$
Define $\rho : S \to T$ by $s \rho = (i, s P_{ji}^{-1}, j)$. It is a
routine exercise to verify that $\rho$ is an isomorphism.

By \cite[Main Theorem]{Ayik99}, $S$ is finitely generated, and so
$T$ is finitely generated. Let
$\sigma : X^+ \to S$
be a choice of semigroup generators for $S$; it follows easily that
$\sigma \rho : X^+ \to T$ is a choice of semigroup generators for
$T$, and that $L_{\sigma}(S) = L_{\sigma \rho}(T)$. Choose now a finite
choice of generators $\tau : Y^+ \to M$ such that $X \subseteq Y$
and $\sigma \rho$ is the restriction of $\tau$ to $X^+$.

Since $P_{ji}$ is a unit, we have $S P_{ji} = S$, so that
$S P_{j' i} \subseteq S P_{ji}$ for every $j' \in J$ with
$P_{j' i} \neq 0$. Now from Theorem~\ref{thm_reestosemi} we
deduce that
$$L_{\sigma}(S) \ = \ L_{\sigma \rho}(T) \ = \ L_\tau(M) \cap \hat{X}^+$$
as required.
\end{proof}

\section{Completely Zero-simple Semigroups}\label{sec_czeros}

Theorems~\ref{thm_semitoreeszero} and~\ref{thm_reestosubsemi}
are of particular interest in one special case. Recall that a finitely
generated semigroup with zero is called \textit{completely zero-simple}
if it has finitely many idempotents and no non-zero ideals. A theorem
of Rees \cite{Rees40} says that every completely zero-simple semigroup
is isomorphic to a Rees matrix semigroup with zero over any of its
non-zero maximal
subgroups, with a sandwich matrix in which every row and every column
contains a non-zero entry. Hence, we obtain the following description
of the loop problem for completely zero-simple semigroups.

\begin{theorem}\label{thm_czeros}
Let $\FamilyF$ be a family of languages closed under union, rational
transduction, product and Kleene closure. Let $S$ be a finitely generated
completely zero-simple semigroup. Then the following are equivalent:
\begin{itemize}
\item $S$ has loop problem in $\FamilyF$;
\item any non-zero maximal subgroup of $S$ has loop problem in $\FamilyF$;
\item every non-zero maximal subgroup of $S$ has loop problem in $\FamilyF$.
\end{itemize}
\end{theorem}
(In fact it is possible to eliminate from the hypothesis of
Theorem~\ref{thm_czeros} the requirements that $\FamilyF$ be closed
under union and product. This can be proved directly using a modification of
the proof of \cite[Theorem~5.4]{KambitesLoop}; the technique makes use of
the group structure of the base semigroup, and so does not seem to lead
to a corresponding strengthening of Theorem~\ref{thm_semitoreeszero}.
Since the proof of the stronger result is rather lengthy, and since
very few interesting language classes are closed under Kleene closure
and transductions but not under product and union, we content ourselves here
with the slightly weaker statement which can be derived easily from
Theorem~\ref{thm_semitoreeszero}.)

In \cite{KambitesLoop} we posed the question of which semigroups have
context-free loop problem. A theorem of Muller and Schupp
\cite{Muller83}, combined with a subsequent result of Dunwoody
\cite{Dunwoody85}, says that a finitely generated group has
context-free word problem exactly if it is virtually free, that is, 
has a free subgroup of finite index. In \cite{KambitesLoop} we applied
this result to give a complete characterization of completely simple
semigroups with context-free loop problem. We are now in a position
to state a corresponding result for completely zero-simple semigroups.
\begin{corollary}
A finitely generated completely zero-simple semigroup has context-free
loop problem if and only if its maximal subgroups are virtually free.
\end{corollary}

\section*{Acknowledgements}

This research was supported by an RCUK Academic Fellowship. The author
would like to thank Kirsty for all her support and encouragement.

\bibliographystyle{plain}

\def\cprime{$'$} \def\cprime{$'$}

\bibliography{mark}
\end{document}